

\def\spaces{\space\space\space\space\space\space\space\space\space\space}
\def\spacess{\message{\spaces\spaces\spaces\spaces\spaces\spaces\spaces}}
\spacess
\spacess
\message{Annals of Mathematics Style: Current Version: 1.0. March 25, 1992}
\spacess
\spacess
\catcode`@=11



\font\elevensc=cmcsc10 scaled\magstephalf
\font\tensc=cmcsc10

\font\eightsc=cmcsc10 scaled800

\font\elevenrm=cmr10 scaled \magstephalf
\font\ninerm=cmr9
\font\eightrm=cmr8
\font\sixrm=cmr6
\font\fiverm=cmr5

\font\eleveni=cmmi10 scaled\magstephalf
\font\ninei=cmmi9
\font\eighti=cmmi8
\font\sixi=cmmi6
\font\fivei=cmmi5
\skewchar\ninei='177 \skewchar\eighti='177 \skewchar\sixi='177
\skewchar\eleveni= '177

\font\elevensy=cmsy10 scaled\magstephalf
\font\ninesy=cmsy9
\font\eightsy=cmsy8
\font\sixsy=cmsy6
\font\fivesy=cmsy5
\skewchar\ninesy='60 \skewchar\eightsy='60 \skewchar\sixsy='60
\skewchar\elevensy='60

\font\eighteenbf=cmbx10 scaled\magstep3

\font\twelvebf=cmbx10 scaled \magstep1
\font\elevenbf=cmbx10 scaled \magstephalf
\font\tenbf=cmbx10
\font\ninebf=cmbx9
\font\eightbf=cmbx8
\font\sixbf=cmbx6
\font\fivebf=cmbx5

\font\elevenit=cmti10 scaled\magstephalf
\font\nineit=cmti9
\font\eightit=cmti8

\font\eighteenmib=cmmi10 scaled \magstep3
\font\twelvemib=cmmib10 scaled \magstep1
\font\elevenmib=cmmib10 scaled\magstephalf
\font\tenmib=cmmib10
\font\eightmib=cmmib10 scaled 800 
\font\sixmib=cmmib10 scaled 600

\font\eighteensyb=cmsy10 scaled \magstep3
\font\twelvesyb=cmbsy10 scaled \magstep1
\font\elevensyb=cmbsy10 scaled \magstephalf
\font\tensyb=cmbsy10 
\font\eightsyb=cmbsy10 scaled 800
\font\sixsyb=cmbsy10 scaled 600
 
\font\elevenex=cmex10 scaled \magstephalf
\font\tenex=cmex10     
\font\eighteenex=cmex10 scaled \magstep3


\def\elevenpoint{\def\rm{\fam0\elevenrm}%
  \textfont0=\elevenrm \scriptfont0=\eightrm \scriptscriptfont0=\sixrm
  \textfont1=\eleveni \scriptfont1=\eighti \scriptscriptfont1=\sixi
  \textfont2=\elevensy \scriptfont2=\eightsy \scriptscriptfont2=\sixsy
  \textfont3=\elevenex \scriptfont3=\elevenex \scriptscriptfont3=\elevenex
  \def\bf{\fam\bffam\elevenbf}%
  \def\it{\fam\itfam\elevenit}%
  \textfont\bffam=\elevenbf \scriptfont\bffam=\eightbf
   \scriptscriptfont\bffam=\sixbf
\normalbaselineskip=13.9pt
  \setbox\strutbox=\hbox{\vrule height9.5pt depth4.4pt width0pt\relax}%
  \normalbaselines\rm}

\elevenpoint 

\def\ninepoint{\def\rm{\fam0\ninerm}%
  \textfont0=\ninerm \scriptfont0=\sixrm \scriptscriptfont0=\fiverm
  \textfont1=\ninei \scriptfont1=\sixi \scriptscriptfont1=\fivei
  \textfont2=\ninesy \scriptfont2=\sixsy \scriptscriptfont2=\fivesy
  \textfont3=\tenex \scriptfont3=\tenex \scriptscriptfont3=\tenex
  \def\it{\fam\itfam\nineit}%
  \textfont\itfam=\nineit
  \def\bf{\fam\bffam\ninebf}%
  \textfont\bffam=\ninebf \scriptfont\bffam=\sixbf
   \scriptscriptfont\bffam=\fivebf
\normalbaselineskip=11pt
  \setbox\strutbox=\hbox{\vrule height8pt depth3pt width0pt\relax}%
  \normalbaselines\rm}

\def\eightpoint{\def\rm{\fam0\eightrm}%
  \textfont0=\eightrm \scriptfont0=\sixrm \scriptscriptfont0=\fiverm
  \textfont1=\eighti \scriptfont1=\sixi \scriptscriptfont1=\fivei
  \textfont2=\eightsy \scriptfont2=\sixsy \scriptscriptfont2=\fivesy
  \textfont3=\tenex \scriptfont3=\tenex \scriptscriptfont3=\tenex
  \def\it{\fam\itfam\eightit}%
  \textfont\itfam=\eightit
  \def\bf{\fam\bffam\eightbf}%
  \textfont\bffam=\eightbf \scriptfont\bffam=\sixbf
   \scriptscriptfont\bffam=\fivebf
\normalbaselineskip=12pt
  \setbox\strutbox=\hbox{\vrule height8.5pt depth3.5pt width0pt\relax}%
  \normalbaselines\rm}


\def\eighteenbold{\def\rm{\fam0\eighteenrm}%
  \textfont0=\eighteenbf \scriptfont0=\twelvebf \scriptscriptfont0=\tenbf
  \textfont1=\eighteenmib \scriptfont1=\twelvemib\scriptscriptfont1=\tenmib
  \textfont2=\eighteensyb \scriptfont2=\twelvesyb\scriptscriptfont2=\tensyb
  \textfont3=\eighteenex \scriptfont3=\tenex \scriptscriptfont3=\tenex
  \def\bf{\fam\bffam\eighteenbf}%
  \textfont\bffam=\eighteenbf \scriptfont\bffam=\twelvebf
   \scriptscriptfont\bffam=\tenbf
\normalbaselineskip=20pt
  \setbox\strutbox=\hbox{\vrule height13.5pt depth6.5pt width0pt\relax}%
\everymath {\fam0 }
\everydisplay {\fam0 }
  \normalbaselines\bf}

\def\elevenbold{\def\rm{\fam0\elevenrm}%
  \textfont0=\elevenbf \scriptfont0=\eightbf \scriptscriptfont0=\sixbf%
  \textfont1=\elevenmib \scriptfont1=\eightmib \scriptscriptfont1=\sixmib%
  \textfont2=\elevensyb \scriptfont2=\eightsyb \scriptscriptfont2=\sixsyb%
  \textfont3=\elevenex \scriptfont3=\elevenex \scriptscriptfont3=\elevenex%
  \def\bf{\fam\bffam\elevenbf}%
  \textfont\bffam=\elevenbf \scriptfont\bffam=\eightbf%
   \scriptscriptfont\bffam=\sixbf%
\normalbaselineskip=14pt%
  \setbox\strutbox=\hbox{\vrule height10pt depth4pt width0pt\relax}%
\everymath {\fam0 }%
\everydisplay {\fam0 }%
  \normalbaselines\bf}


\hsize=31pc
\vsize=48pc

\parindent=22pt
\parskip=0pt

\widowpenalty=10000
\clubpenalty=10000

\topskip=12pt

\skip\footins=20pt
\dimen\footins=3in 

\abovedisplayskip=6.95pt plus3.5pt minus 3pt
\belowdisplayskip=\abovedisplayskip


\voffset=7pt
\def\bindingoffset{\ifodd\pageno \hoffset=-18pt\else\hoffset9pc\fi}

\newif\iftitle

\def\amheadline{\iftitle
\hbox to\hsize{\hss\currannalsline\hss}\else\line{\ifodd\pageno
\hfill\thetitle\hfill\llap{\elevenrm\folio}\else\rlap{\elevenrm\folio}
\hfill\theauthors\hfill\fi}\fi}

\headline={\amheadline}
\footline={\global\titlefalse}
\output={\bindingoffset\plainoutput}


\def\annalsline#1#2{\vfill\eject
\ifodd\pageno\else 
\line{\hfill}
\vfill\eject\fi
\global\titletrue
\def\currannalsline{\eightrm Annals of Mathematics,
{\eightbf#1} (#2), \thepages}}

\def\titleheadline#1{\def\one{#1}\ifx\one\empty\else
\gdef\thetitle{{\frenchspacing%
\let\\ \relax\eightsc\uppercase{#1}}}\fi}

\newif\ifshort

\let\shorttitle\titleheadline

\def\onpages#1#2{\def\thepages{#1--#2}}

\def\thismuchskip[#1]{\vskip#1pt}
\def\ilook{\ifx\next[ \let\go\thismuchskip\else
\let\go\relax\vskip1pt\fi\go}

\def\institution#1{\def\theinstitutions{\vbox{\baselineskip10pt
\def\and{{\eightrm and }}
\def\\{\futurelet\next\ilook}\eightsc #1}}}
\let\institutions\institution

\newwrite\auxfile

\def\startingpage#1{\def\one{#1}\ifx\one\empty\global\pageno=1\else
\global\pageno=#1\fi
\theoremcount=0 \eqcount=0 \sectioncount=0 
\openin1 \jobname.aux \ifeof1 
\onpages{#1}{???}
\else\closein1 \relax\input \jobname.aux
\onpages{#1}{\lastpage}
\fi\immediate\openout\auxfile=\jobname.aux
}

\def\endarticle{\ifRefsUsed\global\RefsUsedfalse%
\else\vskip21pt\theinstitutions%
\nobreak\vskip8pt\vbox{\thereceived\therevised}\fi%
\write\auxfile{\string\def\string\lastpage{\the\pageno}}}

\outer\def\bye{\endarticle\par \vfill \supereject \end}

\newif\ifacks
\long\def\acknowledgements#1{\def\one{#1}\ifx\one\empty\else
\global\ackstrue\vskip-\baselineskip%
\footnote{\ \unskip}{*#1}\fi}

\def\received#1{\def\thereceived{\medbreak
\centerline{\ninerm (Received #1)}}}

\def\title#1{\vbox to76pt{\vfill
\baselineskip=18pt
\parindent=0pt
\overfullrule=0pt
\hyphenpenalty=10000
\everypar={\hskip\parfillskip\relax}
\hbadness=10000
\def\\ {\vskip1sp}
\eighteenbold#1\vskip1sp}
\titleheadline{#1}}

\newif\ifauthor
\def\author#1{\vskip15pt
\hbox to\hsize{\hss\tenrm By \tensc#1\ifacks\global\acksfalse*\fi\hss}
\ifshort\else\xdef\theauthors{{\eightsc\uppercase{#1}}}\fi%
\vskip7pt
\vskip\baselineskip
\global\authortrue\everypar={\global\authorfalse\everypar={}}}

\def\twoauthors#1#2{\vskip15pt
\hbox to\hsize{\hss%
\tenrm By \tensc#1 {\tenrm and} #2\ifacks\global\acksfalse*\fi\hss}
\ifshort\else\xdef\theauthors{{\eightsc\uppercase{#1 and #2}}}\fi
\vskip7pt
\vskip\baselineskip
\global\authortrue\everypar={\global\authorfalse\everypar={}}}


\newcount\theoremcount
\newcount\sectioncount
\newcount\eqcount

\newif\ifspecialnumon

\def\eqnumber=#1 {\global\eqcount=#1 \global\advance\eqcount by-1\relax}
\def\sectionnumber=#1 {\global\sectioncount=#1 
\global\advance\sectioncount by-1\relax}
\def\proclaimnumber=#1 {\global\theoremcount=#1 
\global\advance\theoremcount by-1\relax}

\newif\ifsection
\newif\ifsubsection

\def\intro{\centerline{\bf Introduction}\global\everypar={}
\global\authorfalse\vskip6pt}

\def\elevenboldmath#1{$#1$\egroup}
\def\mathbold{\hbox\bgroup\elevenbold\elevenboldmath}

\def\section#1{\global\theoremcount=0
\global\eqcount=0
\ifauthor\global\authorfalse\else%
\vskip18pt plus 18pt minus 6pt\fi%
{\parindent=0pt\everypar={\hskip\parfillskip}%
\def\\ {\vskip1sp}\elevenpoint\bf%
\ifspecialnumon\global\specialnumonfalse$\rm\spnum$%
\gdef\sectnum{$\rm\spnum$}%
\else\interlinepenalty=10000%
\global\advance\sectioncount by1\relax\the\sectioncount%
\gdef\sectnum{\the\sectioncount}%
\fi.\hskip6pt#1\vrule width 0pt depth12pt}\hskip\parfillskip\break%
\global\sectiontrue%
\everypar={\global\sectionfalse\global\interlinepenalty=0\everypar={}}%
\ignorespaces}


\newif\ifspequation

\let\eqno\leqno 

\newif\ifineqalignno
\let\saveleqalignno\leqalignno
\def\leqalignno{\let\eqnu\Eeqnu\saveleqalignno}

\def\sectandeqnum{%
\ifspecialnumon\global\specialnumonfalse
$\rm\spnum$\gdef\eqnum{$\rm\spnum$}\else\global\firstlettertrue
\global\advance\eqcount by1 
\ifappend\applett\else\the\sectioncount\fi.%
\the\eqcount
\xdef\eqnum{\ifappend\applett\else\the\sectioncount\fi.\the\eqcount}\fi}

\def\eqnu{\leqno{\hbox{\elevenrm\ifspequation\else(\fi\sectandeqnum
\ifspequation\global\spequationfalse\else)\fi}}}

\def\Speqnu{\global\setbox\leqnobox=\hbox{\elevenrm
\ifspequation\else%
(\fi\sectandeqnum\ifspequation\global\spequationfalse\else)\fi}}

\def\Eeqnu{\hbox{\elevenrm
\ifspequation\else%
(\fi\sectandeqnum\ifspequation\global\spequationfalse\else)\fi}}

\newif\iffirstletter
\global\firstlettertrue
\def\eqletter#1{\global\specialnumontrue\iffirstletter\global\firstletterfalse
\global\advance\eqcount by1\fi
\gdef\spnum{\ifappend\applett\else\the\sectioncount\fi.\the\eqcount#1}\eqnu}

\newbox\leqnobox
\def\outsideeqnu#1{\global\setbox\leqnobox=\hbox{#1}}

\def\eatone#1{}

\def\dosplit#1#2{\vskip-.5\abovedisplayskip
\setbox0=\hbox{$\let\eqno\outsideeqnu%
\let\eqnu\Speqnu\let\leqno\outsideeqnu#2$}%
\setbox1\vbox{\noindent\hskip\wd\leqnobox\ifdim\wd\leqnobox>0pt\hskip1em\fi%
$\displaystyle#1\mathstrut$\hskip0pt plus1fill\relax
\vskip1pt
\line{\hfill$\let\eqnu\eatone\let\leqno\eatone%
\displaystyle#2\mathstrut$\ifmathqed~~\qed\fi}}%
\copy1
\ifvoid\leqnobox
\else\dimen0=\ht1 \advance\dimen0 by\dp1
\vskip-\dimen0
\vbox to\dimen0{\vfill
\hbox{\unhbox\leqnobox}
\vfill}
\fi}



\def\testforbreak#1\splitmath #2@{\def\mathtwo{#2}\ifx\mathtwo\empty%
#1$$%
\ifmathqed\vskip-\belowdisplayskip
\setbox0=\vbox{\let\eqno\relax\let\eqnu\relax$\displaystyle#1$}%
\vskip-\ht0\vskip-3.5pt\hbox to\hsize{\hfill\qed}
\vskip\ht0\vskip3.5pt\fi
\else$$\vskip-\belowdisplayskip
\vbox{\dosplit{#1}{\let\eqno\eatone
\let\splitmath\relax#2}}%
\nobreak\vskip.5\belowdisplayskip
\noindent\ignorespaces\fi}


\newif\ifmathqed



\newcount\linenum
\newcount\colnum
\def\spline{\omit&\multispan{\the\colnum}{\hrulefill}\cr}
\def\colcounter{\ifnum\linenum=1\global\advance\colnum by1\fi}

\def\everyline{\noalign{\global\advance\linenum by1\relax}%
\ifnum\linenum=2\spline\fi}

\def\mtable{\bgroup\offinterlineskip
\everycr={\everyline}\global\linenum=0
\halign\bgroup\vrule height 10pt depth 4pt width0pt
\hfill$##$\hfill\hskip6pt\ifnum\linenum>1
\vrule\fi&&\colcounter\hskip12pt\hfill$##$\hfill\hskip12pt\cr}

\def\endmtable{\crcr\egroup\egroup}




\def\xast{*}
\newcount\intable
\newcount\mathcol
\newcount\savemathcol
\newcount\topmathcol
\newdimen\arrayhspace
\newdimen\arrayvspace

\arrayhspace=8pt 
\arrayvspace=12pt 

\newif\ifdollaron

\def\mathalign#1{\def\arg{#1}\ifx\arg\xast%
\let\go\relax\else\let\go\mathalign%
\global\advance\mathcol by1 %
\global\advance\topmathcol by1 %
\expandafter\def\csname  mathcol\the\mathcol\endcsname{#1}%
\fi\go}

\def\arraypickapart#1]#2*{\if#1c \ifmmode\vcenter\else
\global\dollarontrue$\vcenter\fi\else%
\if#1t\vtop\else\if#1b\vbox\fi\fi\fi\bgroup%
\def\one{#2}}

\def\arraystrut{\vrule height .7\arrayvspace depth .3\arrayvspace width 0pt}

\def\array#1#2*{\def\firstarg{#1}%
\if\firstarg[ \def\two{#2} \expandafter\arraypickapart\two*\else%
\ifmmode\vcenter\else\vbox\fi\bgroup \def\one{#1#2}\fi%
\global\everycr={\noalign{\global\mathcol=\savemathcol\relax}}%
\def\\ {\cr}%
\global\advance\intable by1 %
\ifnum\intable=1 \global\mathcol=0 \savemathcol=0 %
\else \global\advance\mathcol by1 \savemathcol=\mathcol\fi%
\expandafter\mathalign\one*%
\mathcol=\savemathcol %
\halign\bgroup&\hskip.5\arrayhspace\arraystrut%
\global\advance\mathcol by1 \relax%
\expandafter\if\csname mathcol\the\mathcol\endcsname r\hfill\else%
\expandafter\if\csname mathcol\the\mathcol\endcsname c\hfill\fi\fi%
$\displaystyle##$%
\expandafter\if\csname mathcol\the\mathcol\endcsname r\else\hfill\fi\relax%
\hskip.5\arrayhspace\cr}

\def\endarray{\crcr\egroup\egroup%
\global\mathcol=\savemathcol %
\global\advance\intable by -1\relax%
\ifnum\intable=0 %
\ifdollaron\global\dollaronfalse $\fi
\loop\ifnum\topmathcol>0 %
\expandafter\def\csname  mathcol\the\topmathcol\endcsname{}%
\global\advance\topmathcol by-1 \repeat%
\global\everycr={}\fi%
}

\def\big#1{{\hbox{$\left#1\vbox to 10pt{}\right.\n@space$}}}
\def\Big#1{{\hbox{$\left#1\vbox to 13pt{}\right.\n@space$}}}
\def\bigg#1{{\hbox{$\left#1\vbox to 16pt{}\right.\n@space$}}}
\def\Bigg#1{{\hbox{$\left#1\vbox to 19pt{}\right.\n@space$}}}


\def\figcaption#1#2#3{\topinsert
\vskip4pt 
\vbox to#3{\vfill}\vskip1sp
\setbox0=\hbox{\eightsc Figure #1.\hskip12pt\eightpoint\rm #2}
\ifdim\wd0>\hsize
\noindent\eightsc Figure #1.\hskip12pt\eightpoint\rm #2\else
\centerline{\eightsc Figure #1.\hskip12pt\eightpoint\rm #2}\fi
\vskip16pt
\endinsert}

\def\wfig#1#2#3{\topinsert
\vskip4pt 
\hbox to\hsize{\hss\vbox{\hrule height .25pt width #3
\hbox to #3{\vrule width .25pt height #2\hfill\vrule width .25pt height #2}
\hrule height.25pt}\hss}
\vskip1sp
\centerline{\eightsc Figure #1}
\vskip16pt
\endinsert}

\def\wfigcaption#1#2#3#4{\topinsert
\vskip4pt 
\hbox to\hsize{\hss\vbox{\hrule height .25pt width #4
\hbox to #4{\vrule width .25pt height #3\hfill\vrule width .25pt height #3}
\hrule height.25pt}\hss}
\vskip1sp
\setbox0=\hbox{\eightsc Figure #1.\hskip12pt\eightpoint\rm #2}
\ifdim\wd0>\hsize
\noindent\eightsc Figure #1.\hskip12pt\eightpoint\rm #2\else
\centerline{\eightsc Figure #1.\hskip12pt\eightpoint\rm #2}\fi
\vskip16pt
\endinsert}

\def\tabcaption#1#2{\vskip6pt
\setbox0=\hbox{\eightsc Table #1.\hskip12pt #2}
\ifdim\wd0>\hsize
\noindent\eightsc Table #1.\hskip12pt\eightpoint #2
\else
\centerline{\eightsc Table #1.\hskip12pt\eightpoint #2}
\fi
\vskip6pt}

\def\endinsert{\egroup\if@mid\dimen@\ht\z@\advance\dimen@\dp\z@ 
\advance\dimen@ 12\p@\advance\dimen@\pagetotal\ifdim\dimen@ >\pagegoal 
\@midfalse\p@gefalse\fi\fi\if@mid\smallskip\box\z@\bigbreak\else
\insert\topins{\penalty 100 \splittopskip\z@skip\splitmaxdepth\maxdimen
\floatingpenalty\z@\ifp@ge\dimen@\dp\z@\vbox to\vsize {\unvbox \z@ 
\kern -\dimen@ }\else\box\z@\nobreak\smallskip\fi}\fi\endgroup}

\def\pagecontents{
\ifvoid\topins \else\iftitle\else 
\unvbox \topins \fi\fi \dimen@ =\dp \@cclv \unvbox 
\@cclv 
\ifvoid\topins\else\iftitle\unvbox\topins\fi\fi
\ifvoid \footins \else \vskip \skip \footins \footnoterule 
\unvbox \footins \fi \ifr@ggedbottom \kern -\dimen@ \vfil \fi}


\newif\ifappend

\def\appendix#1#2{\def\applett{#1}\def\two{#2}%
\global\appendtrue
\global\theoremcount=0
\global\eqcount=0
\vskip18pt plus 18pt
\vbox{\parindent=0pt
\everypar={\hskip\parfillskip}
\def\\ {\vskip1sp}\elevenpoint\bf Appendix%
\ifx\applett\empty\gdef\applett{A}\ifx\two\empty\else.\fi%
\else\ #1.\fi\hskip6pt#2\vskip12pt}%
\global\sectiontrue%
\everypar={\global\sectionfalse\everypar={}}\nobreak\ignorespaces}

\newif\ifRefsUsed
\long\def\references{\global\RefsUsedtrue\vskip21pt
\theinstitutions
\everypar={}\global\bibnum=0
\vskip20pt\goodbreak\bgroup
\vbox{\centerline{\eightsc References}\vskip6pt}%
\ifdim\maxbibwidth>0pt\leftskip=\maxbibwidth\else\leftskip=18pt\fi%
\ifdim\maxbibwidth>0pt\parindent=-\maxbibwidth\else\parindent=-18pt\fi%
\everypar={\amref}%
\ninepoint
\frenchspacing
\nobreak\ignorespaces}

\def\endreferences{\vskip1sp\egroup\everypar={}%
\nobreak\vskip8pt\vbox{\thereceived\therevised}}

\newcount\bibnum

\def\amref#1 {\def\one{#1}\global\advance\bibnum by1%
\immediate\write\auxfile{\string\expandafter\string\def\string\csname
\space #1\string\endcsname{[\the\bibnum]}}%
\leavevmode\hbox to18pt{\hbox to13.2pt{\hss[\the\bibnum]}\hfill}}

\def\bibline{\hbox to30pt{\hrulefill}\/\/}

\def\name#1{{\eightsc#1}}

\newdimen\maxbibwidth
\def\AuthorRefNames [#1] {\def\amref{\spamref}
\setbox0=\hbox{[#1] }\global\maxbibwidth=\wd0\relax}

\def\spamref[#1] {\leavevmode\hbox to\maxbibwidth{\hss[#1]\hfill}}


\def\footnoterule{\kern-3pt\hrule width1in height.5pt\kern2.5pt}

\let\savefootnote\footnote
\def\footnote#1#2{%
\savefootnote{#1}{{\eightpoint\normalbaselineskip11pt
\normalbaselines#2}}}

\def\vfootnote#1{%
\insert \footins \bgroup \eightpoint\normalbaselineskip11pt\normalbaselines
\interlinepenalty \interfootnotelinepenalty
\splittopskip \ht \strutbox \splitmaxdepth \dp \strutbox \floatingpenalty 
\@MM \leftskip \z@skip \rightskip \z@skip \spaceskip \z@skip 
\xspaceskip \z@skip
{#1}$\,$\footstrut \futurelet \next \fo@t}


\newif\iffirstadded
\newif\ifadded

\def\addedlett{}

\def\alltheoremnums{%
\ifspecialnumon\global\specialnumonfalse
\ifadded\global\addedfalse
\iffirstadded\global\firstaddedfalse
\global\advance\theoremcount by1 \fi
\ifappend\applett\else\the\sectioncount\fi.\the\theoremcount\addedlett%
\xdef\theoremnum{\ifappend\applett\else\the\sectioncount\fi.%
\the\theoremcount\addedlett}%
\else$\rm\spnum$\def\theoremnum{$\rm\spnum$}\fi%
\else\global\firstaddedtrue
\global\advance\theoremcount by1 
\ifappend\applett\else\the\sectioncount\fi.\the\theoremcount%
\xdef\theoremnum{\ifappend\applett\else\the\sectioncount\fi.%
\the\theoremcount}\fi}

\def\allcorolnums{%
\ifspecialnumon\global\specialnumonfalse
\ifadded\global\addedfalse
\iffirstadded\global\firstaddedfalse
\global\advance\corolcount by1 \fi
\the\corolcount\addedlett%
\else$\rm\spnum$\def\corolnum{$\rm\spnum$}\fi%
\else\global\advance\corolcount by1 
\the\corolcount\fi}


\newcount\corolcount
\def\xcorol{Corollary}
\def\xtheorem{Theorem}
\def\xmaintheorem{Main Theorem}

\newif\ifthtitle

\let\saverparen)
\let\savelparen(
\def\rmparenl{{\rm(}}
\def\rmparenr{{\rm\/)}}
{
\catcode`(=13
\catcode`)=13
\gdef\makeparensRM{
\catcode`(=13
\catcode`)=13
\let(=\rmparenl
\let)=\rmparenr
\everymath{\let(\savelparen
\let)\saverparen}
\everydisplay{\let(\savelparen
\let)\saverparen}}
}

\medskipamount=8pt plus.1\baselineskip minus.05\baselineskip
\def\proclaim#1{\vskip-\lastskip
\def\one{#1}\ifx\one\xtheorem\global\corolcount=0\fi
\ifsection\global\sectionfalse\vskip-6pt\fi
\medskip
{\elevensc#1}%
\ifx\one\xmaintheorem\global\corolcount=0
\gdef\theoremnum{Main Theorem}\else%
\ifx\one\xcorol\ \allcorolnums\else\ \alltheoremnums\fi\fi%
\ifthtitle\ \global\thtitlefalse{\rm(\thethtitle)}\fi.%
\hskip.5pc\bgroup\makeparensRM
\it\ignorespaces}

\def\nonumproclaim#1{\vskip-\lastskip
\def\one{#1}\ifx\one\xtheorem\global\corolcount=0\fi
\ifsection\global\sectionfalse\vskip-6pt\fi
\medskip
{\elevensc#1}.\ifx\one\xmaintheorem\global\corolcount=0
\gdef\theoremnum{Main Theorem}\fi\hskip.5pc%
\bgroup\makeparensRM\it\ignorespaces}

\def\endproclaim{\egroup\medskip}


\def\xproof{Proof}
\def\xremark{Remark}
\def\xcase{Case}
\def\xsubcase{Subcase}
\def\xconjecture{Conjecture}
\def\xstep{Step}
\def\xof{of}

\def\deconstruct#1 #2 #3 #4 #5 @{\def\one{#1}\def\two{#2}\def\three{#3}%
\def\four{#4}%
\ifx\two\empty #1\else%
\ifx\one\xproof%
\ifx\two\xof%
  \ifx\three\xcorol Proof of Corollary \rm#4\else%
     \ifx\three\xtheorem Proof of Theorem \rm#4\else\xone\fi%
  \fi\fi%
\else\xone\fi\fi.}

\def\pickup#1 {\def\this{#1}%
\ifx\this\xproof\global\let\go\demoproof
\global\let\enddemo\endproof\else
\ifx\this\xremark\global\let\go\demoremark\else
\ifx\this\xcase\global\let\go\demostep\else
\ifx\this\xsubcase\global\let\go\demostep\else
\ifx\this\xconjecture\global\let\go\demostep\else
\ifx\this\xstep\global\let\go\demostep\else
\global\let\go\demoproof\fi\fi\fi\fi\fi\fi}

\def\demo#1{\vskip-\lastskip
\ifsection\global\sectionfalse\vskip-6pt\fi
\def\one{#1 }\def\two{#1*}%
\setbox0=\hbox{\expandafter\pickup\one}\expandafter\go\two}

\def\numbereddemo#1{\vskip-\lastskip
\ifsection\global\sectionfalse\vskip-6pt\fi
\def\two{#1*}%
\expandafter\demoremark\two}

\def\demoproof#1*{\medskip\def\xone{#1}
{\ignorespaces\it\expandafter\deconstruct\xone {} {} {} {} {} @%
\unskip\hskip6pt}\rm\ignorespaces}

\def\demoremark#1*{\medskip
{\it\ignorespaces#1\/} \alltheoremnums\unskip.\hskip1pc\rm\ignorespaces}

\def\demostep#1 #2*{\vskip4pt
{\it\ignorespaces#1\/} #2.\hskip1pc\rm\ignorespaces}

\def\enddemo{\medskip}

\def\endproof{\ifmathqed\global\mathqedfalse\medskip\else
\parfillskip=0pt~~\hfill\qed\medskip
\fi\global\parfillskip0pt plus 1fil\relax
\gdef\enddemo{\medskip}}

\def\qed{\vbox{\hrule\hbox{\vrule height6pt\hskip6pt\vrule}\hrule}}








\def\stripbs#1#2*{\def\one{#2}}
\def\newdef#1#2#{\expandafter\stripbs\string#1*
\defining{#1}{#2}}

\def\defining#1#2#3{\ifdefined{\one}{\message{<<<Sorry, 
\string#1\space
is already defined. Please use a new name.>>>}}
{\expandafter\gdef#1#2{#3}}}

\def\emptyspace{ }
\def\nextthing{}
\def\newline{***}
\def\eatone#1{ }

\def\lookatspace#1{\ifcat\noexpand#1\ \else%
\gdef\nextthing{}\xdef\next{#1}%
\ifx\next\emptyspace%
\let\nextthing\emptyspace\else\ifx\next\newline%
\gdef\nextthing{\eatone}\fi\fi\fi\egroup\nextthing#1}

{\catcode`\^^M=\active%
\gdef\spacer{\bgroup\catcode`\^^M=\active%
\let^^M=\newline\obeyspaces\lookatspace}}

\def\ref#1{\seeifdefined{#1}#1\spacer}


\def\seeifdefined#1{\expandafter\stripbs\string#1*%
\crorefdefining{#1}}

\newif\ifcromessage
\global\cromessagetrue

\def\crorefdefining#1{\ifdefined{\one}{}
{\ifcromessage\global\cromessagefalse%
\message{\spaces\spaces\spaces\spaces\spaces\spaces\spaces}%
\message{<Undefined reference.}%
\message{Please TeX file once more to have accurate cross-references.>}%
\message{\spaces\spaces\spaces\spaces\spaces\spaces\spaces}\fi%
\expandafter\gdef#1{??}}}

\def\label#1#2*{\gdef\ctest{#2}\seeifdefined{#1}%
\ifx\empty\ctest\write\auxfile{\noexpand\def\noexpand#1{\tocstuff}}%
\else\def\ctemp{#2}%
\immediate\write\auxfile{\noexpand\def\noexpand#1{\ctemp}}\fi}

\def\ifdefined#1#2#3{\expandafter\ifx\csname#1\endcsname\relax#3\else#2\fi}




\def\articlecontents{
\vskip20pt\centerline{\bf Table of Contents}\everypar={}\vskip6pt
\bgroup \leftskip=3pc \parindent=-2pc 
\def\item##1{\vskip1sp\indent\hbox to2pc{##1.\hfill}}}

\def\endcontents{\vskip1sp\leftskip=0pt\egroup}

\def\journalcontents{\vfill\eject
\def\currannalsline{\hfill}
\global\titletrue
\vglue3.5pc
\centerline{\tensc\hskip12pt TABLE OF CONTENTS}\everypar={}\vskip30pt
\bgroup \leftskip=34pt \rightskip=-12pt \parindent=-22pt 
  \def\\ {\vskip1sp\noindent}
\def\pagenum##1{\unskip\parfillskip=0pt\dotfill##1\vskip1sp
\parfillskip=0pt plus 1fil\relax}
\def\name##1{{\tensc##1}}}

\def\indexcontents{\vfill\eject
\leftline{\tensc \quad Annals of Mathematics}
\vglue-2truept
\leftline{\quad\phantom{AN}\ninerm   Vol.\ 148, Nos.\ 1, 2, 3}
\vglue-2truept
\leftline{\quad\phantom{ANNALS} \ninerm 1998}
\def\currannalsline{\hfill} 
\def\amheadline{\iftitle
\hbox to\hsize{\hss\currannalsline\hss}\else\line{\ifodd\pageno
\hfill\thetitle\hfill\llap{\elevenrm\folio}\else\rlap{\elevenrm\folio}
\hfill\fi}\fi}

\headline={\amheadline}
\global\titletrue
\vglue2.5pc 
\centerline{\tensc\hskip12pt INDEX}\everypar={}\vskip18pt
\bgroup \leftskip=34pt \rightskip=-12pt \parindent=-22pt 
  \def\\ {\vskip1sp\noindent}
\def\pagenum##1{\unskip\parfillskip=0pt\dotfill##1\vskip1sp
\parfillskip=0pt plus 1fil\relax}
\def\name##1{{\tensc##1}}}


\institution{}
\onpages{0}{0}
\def\lastpage{???}
\def\thetitle{Title ???}
\def\theauthors{Authors ???}
\def\thereceived{}
\def\therevised{}

\catcode`\@=12

\annalsline{151}{2000}
\received{April 21, 1999}
\startingpage{375}

\font\emi= cmmi10 scaled 1700 
\font\eightrm=cmr10  
\font\sixrm=cmr6
\def\Z{{\bf Z}}
\def\Q{{\bf Q}}

\def\N{{\bf N}}
\def\F{{\bf F}}

\def\?{{\cal P}}
\def\B{{\cal B}}
\title{A proof of Pisot's {\emi d}\hskip1pt\raise7pt\hbox{\eightrm th} root conjecture} 
 \def\titleheadline#1{\def\one{#1}\ifx\one\empty\else
\gdef\thetitle{{\frenchspacing%
\let\\ \relax
{#1}}}\fi}
\newif\ifshort

\let\shorttitle\titleheadline
\shorttitle{ \eightsc\uppercase{A proof of Pisot's} {\eightpoint \it d}
\hskip-2pt\raise3pt\hbox{{\sixrm th}}  \eightsc\uppercase{conjecture}}
\def\pmb#1{\setbox0=\hbox{$#1$}%
		\kern-.025em\copy0\kern-\wd0
		\kern.05em\copy0\kern-\wd0
		\kern-.025em\raise.0433em\box0 }
\acknowledgements{}
 \author{Umberto Zannier}
 \institutions{Istituto Universitario di Architettura di Venezia, Venezia, Italy\\
 {\eightpoint {\it E-mail address\/}: zannier@iuav.unive.it}}

\bigbreak 
\centerline{\bf Abstract} 
\bigbreak
Let $\{b(n):n\in\N\}$ be the sequence
of coefficients in the Taylor expansion of a rational function
$R(X)\in\Q(X)$ and
suppose that $b(n)$ is a perfect $d^{\rm th}$ power for all large $n$. A
conjecture of
Pisot states that one can choose a $d^{\rm th}$ root $a(n)$ of $b(n)$ such that $\sum
a(n)X^n$ is also a rational function. Actually, this is the fundamental
case of an
analogous statement formulated for fields more general than $\Q$. A number
of papers
have been devoted to various special cases. In this note we shall
completely settle
the general case.

\phantom{foothurts}
\medbreak
 \intro
\medbreak
 Let $R(X)=\sum_{n=0}^\infty b(n)X^n$ represent a
rational function in $k(X)$, where ${\rm char}\ k=0$. (This is equivalent to a
linear recurrence $\sum_{i=0}^tc_ib(n+i)=0$, valid for large $n$.) Suppose that
there exists a  field $F\supset k$, finitely generated over $\Q$, such that
$b(n)$
is a perfect $d^{\rm th}$ power in $F$ for all large $n$. Then, it is a
generalization of a
conjecture attributed to  Pisot (see e.g. [B], [RvdP] and \S 6.4  of
[vdP1]) that
there exists a sequence $\{a(n)\}$ such that $a(n)^d=b(n)$ for large $n$
and $\tilde
R(X):=\sum a(n)X^n$ is again rational. 

Let $\beta_1,\ldots ,\beta_h$ be the poles of $R(X)$, of multiplicities
$m_1,\ldots
,m_h$ respectively. Then it is well-known that for large $n$, $b(n)$ may be
expressed
as an  exponential polynomial
$$
b(n)=\sum_{i=1}^hB_i(n)\beta_i^n\eqno(1)
$$
for polynomials $B_i\in \bar k[T]$ of degrees resp. $\le m_i-1$. (The
$\beta_i$ are
classically called the {\it roots} of the exponential polynomial.) The conjecture predicts that under the stated
assumptions the right side of
(1) is identically the $d^{\rm th}$ power of a function of the same form. 

Rumely and van der Poorten [RvdP, \S 7] refer to the above statement as the
{\it
Generalized Pisot $d^{\rm th}$ root conjecture}.  In [RvdP, \S 6]  they use
specialization arguments to prove that it is actually sufficient to deal
with the
case when $k$ is a number field. Also, generalizing a previous argument in
[PZ], they
prove  that it is also possible to assume $m_i=1$ for $1\le i\le h$,
namely that all the $B_i$'s are constant (see [RvdP, Prop. 3]). Theorem 2 of
[RvdP] collects these facts. They
prove in Theorem 1 that in the fundamental  number field case, the conjecture is true provided there exists a
unique pole
$\beta_i$ of maximal or minimal absolute value. This
condition, though
rather weak, plays a crucial role at one step in their arguments. In fact, the
same type of condition, leading to the so-called {\it dominant root
method},  had
presented an obstacle also in several related investigations (e.g. on
Pisot conjecture on the {\it Hadamard quotient}, studied by Pisot and finally
solved  by van der Poorten [vdP2] after an incomplete argument by Pourchet
[P].) For
a proof of the $d^{\rm th}$ root conjecture under different additional
assumptions, see
[B]. 

Here we shall completely settle the fundamental case when $k$ is a number
field, by
means of a method  entirely different from those mentioned above. As we
have noticed,
the results in [RvdP] allow us  to assume that the $B_i$ are constant
without loss of
generality. However, since this assumption would only slightly simplify our
proofs,
we shall avoid it. By enlarging
$k$ we may suppose  that $F=k$ and that the condition on $b(n)$ is true for
all $n$.
We have the following 

\nonumproclaim{Theorem} 
\hskip-12pt Let $b(n)=\sum_{i=1}^hB_i(n)\beta_i^n$ be an
exponential
polynomial{\rm ,} where $B_i\in k[T]$ and $\beta_i\in k${\rm ,} for a   field $k$ finitely
generated over $\Q${\rm .} Assume that $b(n)$ is a $d^{\rm th}$ power in $k$ for $n\in\N${\rm .}
Then there exists an exponential polynomial
$a(n)=\sum_{i=1}^rA_i(n)\alpha_i^n$,
$A_i\in\bar k[T]${\rm ,} $\alpha_i\in\bar k${\rm ,} such that $b(n)=a(n)^d$ for all
$n\in\N${\rm .}
\endproclaim

By what we have recalled above, it is sufficient to deal with the case when
$k$ is
a number field, as we shall assume from now  on. We remark that [CZ]
contains proofs
of analogous statements assuming the much weaker fact that {\it $b(n)$ is a
$d^{\rm th}$
power for infinitely many $n\in\N$}. However, that method again requires
(at least)
the existence of a dominant $\beta_i$. On the other hand, the present
method, of
completely different nature,  does not yield results under the weaker
assumption. 

\demo{{R}emark {\rm  1}} The proof will give in fact the (apparently) stronger
result stating that if the conclusion is not true, we may find an
arithmetic progressions ${\cal A}$ and a prime ideal $\pi$ of
$k$ such that  $b(n)$ is not  a $d^{\rm th}$ power modulo $\pi$, for each $n\in{\cal
A}$. In particular, it suffices that the assumption holds for a set of positive
integers intersecting  every arithmetic progression.\enddemo

\demo{{R}emark {\rm  2}} It should be rather straightforward to adapt the
arguments to more general equations $f(y,b(n))=0$, $f\in
k[x,y]$, monic
in~$y$. Assuming that for each $n\in \N$ the equation has a solution in
$k$, one
would obtain that there exists an identical solution which is  an
exponential polynomial. (This would provide in particular a proof of
conjectures
mentioned in [vdP3].) \enddemo

\demo{{P}roof of theorem}  Here is a (very brief) outline of the method.
First, we shall look at what happens upon replacing in the exponential
polynomial the
numbers $\beta_i^n$ with suitable roots of unity  of order dividing $p-1$,
for a
sufficiently large prime $p\equiv 1\pmod d$. Assuming the conclusion to be
false, we
shall use congruences  modulo $p$ to show that the resulting algebraic number
will not be a $d^{\rm th}$ power in the corresponding cyclotomic extension of
$k$. (More
precisely, we shall use the Lang-Weil theorem for the number of points on
varieties
over finite fields.) Now   Cebotarev's theorem will show that, for
infinitely many
prime ideals, the reduction of this algebraic integer will continue not to be a
$d^{\rm th}$ power in the residue field. By taking a suitable $n\in \N$, the
effect of the
reduction will just consist of  replacing  the mentioned roots of unity
with the
$\beta_i^n$. This will contradict the assumption that $b(n)$ is a $d^{\rm th}$
power in
$k$. 

To begin the proof, we point out at once that  it suffices to deal
with the case when $d$ is a prime number and $k$ is normal over $\Q$, as we
shall
assume from now on.

We may assume that no
$\beta_i$ is zero, and we consider the multiplicative subgroup
$\Gamma\subset k^*$
generated by the $\beta_i$'s. We begin by proving that it is sufficient to deal
with the case when $\Gamma$ is torsion-free. Let $N$ be the  order of the
torsion
subgroup of $\Gamma$. For $r=0,1,\ldots ,N-1$ we consider the exponential
polynomial
$b_r(n):=b(r+Nn)$. For an exponential polynomial $c(n)=\sum c_i(n)\gamma_i^n$,
choose in some way $N^{\rm th}$ roots $\tilde\gamma_i$ for the $\gamma_i$'s and put
$\tilde c(n)=\sum c_i({n\over N})\tilde\gamma_i^n$.  We have $\tilde
c(nN)=c(n)$ for
$n\in \N$.

Suppose now that we can
prove the theorem for each $b_r(n)$, so that $b_r(n)=a_r(n)^d$ for suitable
exponential polynomials $a_r(n)$, $0\le r\le N-1$. Define
$\phi(n)=(1/N)\sum_{\theta^N=1}\theta^n$, so $\phi(n)$ is $1$ for $N|n$ and $0$
otherwise. Then we find that
$$
b(n)=\left
(\sum_{r=0}^{N-1}\phi(n-r)\tilde a_r(n-r)\right )^d.
$$
In fact, if $n=s+mN$, for integers $m,s$, $0\le s\le N-1$, the right side
equals $\tilde a_s(mN)^d=a_s(m)^d=b(s+mN)=b(n)$, as required. Since the
expression
into brackets in the right side  is plainly an exponential polynomial, the
theorem is
proved for $b(n)$.

Therefore it suffices to prove the theorem for the $b_r(n)$. On the other
hand,
the roots of $b_r(n)=b(r+nN)$ are just the $\beta_i^N$, which clearly
generate in
$k^*$ the torsion-free group $\Gamma^N$.  Hence the above claim  follows.
Namely, in proving the theorem we can and shall assume that $\Gamma$ is
torsion-free.

We let $\gamma_1,\ldots ,\gamma_r\in k^*$ be multiplicatively independent
(i.e. free)
generators for $\Gamma$. Then we may write
$$
b(n)=f(n,\gamma_1^n,\ldots ,\gamma_r^n),\eqno(2)
$$
for a suitable rational function $f\in k[X_0,X_1,X_1^{-1},\ldots
,X_r,X_r^{-1}]$.
Multiplying  by a  power of $(\gamma_1\cdots\gamma_r)^{nd}$ does not
affect the assumptions nor the conclusion, so we can assume that $f$ is a
polynomial in all the variables.  

At this point we observe a technical fact, which  shall be useful later.
Let $D$ be a
positive integer and consider the polynomial
$$
\eqalign{
\noalign{\vskip6pt}
F(X_0,X_1,\ldots ,X_r)=F_D(X_0,X_1,\ldots ,X_r):=f(X_0,X_1^D,\ldots ,X_r^D).\cr
\noalign{\vskip6pt}
}
$$
Suppose first that $F$ is a $d^{\rm th}$ power in $\overline\Q [X_0,\ldots
,X_r]$, say
$F=G^d$. Observe that $F$ is invariant by every substitution $\phi$,
$X_i\mapsto\theta_i X_i$ where $\theta_i$ are $D^{\rm th}$ roots of unity and
$\theta_0=1$. Therefore we have, for each such $\phi$,
$$\eqalign{
\noalign{\vskip6pt}
\phi(G)=\mu_\phi G,\cr
\noalign{\vskip6pt}
}
$$
where $\mu_\phi^d=1$. From these equations we see at once that the quotient
of any
two monomials in $G$ is invariant by all the substitutions $\phi$. This
means that
$G$ is of the form $X_1^{a_1}\cdots X_r^{a_r}H(X_0,X_1^D,\ldots ,X_r^D)$, for
 $a_1,\ldots ,a_r\in\Z$ and a polynomial $H\in \overline\Q [X_0,\ldots ,X_r]$.
Since $F=G^d$, we have $D|a_id$ for all $i\ge 1$. Therefore  $f(X_0,\ldots
,X_r)$
may be written in  the form $X_1^{b_1}\cdots X_r^{b_r}H^d(X_0,\ldots
,X_r)$, where
$a_id=b_iD$. Now, if $\alpha$ is a $d^{\rm th}$ root of $\gamma_1^{b_1}\cdots
\gamma_r^{b_r}$ in $\bar k$,  equation (2) gives
$b(n)=(\alpha^nH(n,\gamma_1^n,\ldots,\gamma_r^n))^d$, providing a proof of the
theorem.

In conclusion, we may suppose that for each positive integer $D$, $F_D$ is
not a
$d^{\rm th}$ power  in $\overline\Q [X_0,\ldots ,X_r]$. We proceed to
derive a contradiction. 
\bigbreak

Before going on we shall establish the multiplicative independence of suitable
elements of $k^*$, modulo certain  powers. The following lemma will be
helpful; it
appears (in a more general form) in a forthcoming joint paper with
E. Bombieri and D. Masser [BMZ, L. 6] (and perhaps elsewhere), but for
completeness
we give here the short proof of the case we need.\enddemo

\nonumproclaim{Lemma 1}  Let  $\ell$ be a prime number such that $k$
does not
contain a primitive $\ell^{\rm th}$ root of unity{\rm .} Suppose that $u$ lies in an
abelian extension of $k$ and that $u^\ell\in k^*${\rm .}
Then there exists a root of unity $\omega$ such
that $\omega u\in k^*${\rm .}
\endproclaim

\demo{Proof} We have $g(u)=0$, where $g(X):=X^\ell -u^\ell\in k[X]$.
Suppose first that $g$ is reducible over $k$. Then it is classical that
 $u^\ell =v^\ell$ for some
$v\in k^*$. In this case $u=\zeta v$ for some $\ell^{\rm th}$ root of unity $\zeta$,
concluding the proof.

Therefore we may assume that $g$ is irreducible, so $[k(u):k]=\ell$.
Observe that
$k(u)/k$ is Galois, as a subextension of an abelian extension of $k$. Let
$\sigma$ be
a nontrivial element of the Galois group. Since $u^\ell\in k$ we have
$\sigma(u)=\zeta u$ for some nontrivial (whence primitive) $\ell^{\rm th}$ root
of unity
$\zeta$. This shows that $\zeta\in k(u)$. Therefore $[k(\zeta):k]$ divides both
$[k(u):k]=\ell$ and $[\Q(\zeta):\Q]=\ell -1$, so $\zeta\in k$, a contradiction.
\enddemo

We now choose once and for all a large prime number $\beta$, such that $\beta
,\gamma_1,\ldots ,\gamma_r$ are multiplicatively independent. Also, let
$h_1,\ldots ,h_r$ be integers (which shall be specified later) and put
\vglue-2pt
$$
\delta_i:=\gamma_i\beta^{-h_i},\qquad i=1,\ldots ,r.\eqno(3)
$$

\vglue6pt \noindent Define $\Q_c$ to be the maximal cyclotomic extension of $\Q$ and
$k_c:=k\Q_c$. We have 

\nonumproclaim{Lemma 2} \hskip-9pt  There exists a number $L$ depending only on
$k$ and $\beta ,\gamma_1,\ldots ,\gamma_r${\rm ,} with the following property{\rm .}
Let $M$ be a positive integer  and let
$\ell$ be a prime number $>L${\rm .} Then $\beta^M$ does not lie in the group
generated by the $\delta_i$ and by ${k_c^*}^{\ell M}${\rm .}
\endproclaim

\demo{Proof} Assume that the conclusion is false and that we
have integers $a_1,\ldots ,a_r$ such that
$$
\beta^M \delta_1^{a_1}\cdots\delta_r^{a_r}=v^{\ell M},
$$
where $v\in k_c^*$. Suppose also that this $M$ is minimal among such equalities.   We write it  as
$$
\beta^{a_0}\gamma_1^{a_1}\ldots \gamma_r^{a_r}=v^{\ell M},
$$
where $a_0=M-\sum_{i=1}^ra_ih_i$. The left side lies in $k^*$.
Therefore, if $\ell$ is large enough with respect to $k$, we may apply
Lemma 1 and
suppose that $u:=v^M\break =\omega\nu$ for a root of unity $\omega$ and a $\nu\in
k^*$.

Let $Q$ be a positive integer. By a well-known easy result in diophantine
approximation, we may find a positive integer $b\le Q^{r+1}$ and integers
$p_i$ such
that, putting $b_i:=ba_i-\ell p_i$, we have $|b_i|\le \ell Q^{-1}$ for
$i=0,1,\ldots
,r$. (Apply, e.g., Thm.~1A, p. 27 of [Schm1].)

Raising both sides of the last displayed equation to the $b^{\rm th}$ power we obtain
$$
\beta^{b_0}\gamma_1^{b_1}\cdots\gamma_r^{b_r}=\omega^{\ell b}\rho^\ell
$$
with  $\rho\in k^*$ (we take
$\rho =\nu^b(\beta^{p_0}\gamma_1^{p_1}\cdots\gamma_r^{p_r})^{-1}$). Taking Weil
heights we get
$$
\ell h(\rho)\le \max_i(|b_i|)\left
(h(\beta)+\sum_{i=1}^rh(\gamma_i)\right )\ll{\ell\over Q},
$$
where the implied constant depends only on $\beta$ and the $\gamma_i$'s.
Choosing,
e.g., $Q=\lceil{\ell^{1\over r+2}}\rceil$, we see that $h(\rho)\ll
\ell^{-1\over r+2}$
and $0<b\le Q^{r+1}<\ell$ for large $\ell$. On the other hand the Weil
height has a
positive lower bound on $k^*$, outside  the roots of unity. For large
enough $\ell$
our inequality thus forces $\rho$ to be a root of unity. From the
multiplicative
independence of $\beta ,\gamma_1,\ldots ,\gamma_r$ we thus obtain
$b_0=b_1=\ldots
=b_r=0$. Since $1\le b<\ell$ we deduce that  $a_i$ is a multiple of $\ell$
for all
$i=0,1,\ldots ,r$. In particular $\ell$ divides $M$. Writing $M=\ell
M'$, $a_i=\ell a_i'$ we have
$$
\beta^{M'} \delta_1^{a_1'}\cdots\delta_r^{a_r'}=(v\theta)^{\ell M'}
$$
for some root of unity $\theta$. This however contradicts the minimality
of $M$ and concludes the proof. 
\enddemo

In the sequel $\zeta_s$ will denote a primitive $s^{\rm th}$ root of unity.  Take $L$
as in Lemma 2; without loss we can assume $L>d$. Let $D$ be an integer
divisible by
all the primes $\le L$. We may also
assume that $\Q(\zeta_D)$ contains the maximal cyclotomic subfield
$\Q_c\cap k$ of
$k$.

Define $k_1=k(\zeta_D)$. If $\ell\le L$ is a prime we have that
$k_1(\zeta_{\ell
D}) /k_1$ is cyclic of order $\ell$. In fact, $\ell |D$ and
$[k_1(\zeta_{\ell D}):k_1]=[\Q(\zeta_{\ell D}):\Q(\zeta_D)]=\ell$; this is
because
$[\Q(\zeta_{\ell D}):\Q_c\cap k]=[k(\zeta_{\ell D}):k]$ and
$[\Q(\zeta_D):\Q_c\cap
k]=[k_1:k]$.

Let $k_2$ be the compositum of the fields $k_1(\zeta_{\ell D})$ for all primes
$\ell \le L$. Then the Galois group ${\rm Gal}(k_2/k_1)$ is a product of cyclic
groups of
pairwise coprime orders and is thus cyclic. Pick a generator  $\sigma$ of
${\rm Gal}(k_2/k_1)$. Then  $\sigma$ does not induce the identity on any field
$k_1(\zeta_{\ell D})$ as above. 

Let $p$ be a prime number which splits completely in $k_1/\Q$, and such
that, for
some $\B$ in $k_2$ above $p$, the Frobenius automorphism of $\B$ in
${\rm Gal}(k_2/k_1)$
is $\sigma$. Cebotarev's theorem guarantees the existence of infinitely
many such
primes. Since $\Q(\zeta_D)\subset k_1$, $p$ splits completely also in
$\Q(\zeta_D)/\Q$ and hence $p\equiv 1\pmod D$.

On the other hand, $p$ cannot split completely in any $\Q(\zeta_{\ell D})$
for a
prime $\ell \le L$. (Otherwise it would split completely in
$k_1(\zeta_{\ell D})$ and
$\sigma$ would be the identity on this field.) In turn, this means that
$p\not\equiv 1\pmod{\ell D}$. In
conclusion we get infinitely many primes $p$ with the following
properties.\medbreak

\item{(a)} $p=1+Dm$, where the least prime factor of $m$ is $>L>d$.
\smallbreak
\item{(b)} $p$ splits completely in $k_1$.\medbreak

Put, as above, $F(X_0,\ldots ,X_r)=f(X_0,X_1^D,\ldots ,X_r^D)$. We have
seen that it
is possible to assume that $F$ is not a $d^{\rm th}$ power in $\overline\Q[X_0,\ldots
,X_r]$. Therefore (by Kummer theory) the polynomial
$$
S(X_0,\ldots ,X_r,Y)=Y^d-F(X_0,\ldots ,X_n)\eqno(4)
$$
is absolutely irreducible. Pick now a large prime $p$ with the properties
(a),(b) and consider a prime ideal factor ${\cal P}$ of $p$ in $k_1$. For
large $p$
we may look at the reduction $\tilde S\in \F_p[X_0,\ldots ,X_r,Y]$ of $S$
modulo
${\cal P}$. Also, for large $p$ a well-known theorem due to Ostrowski (see e.g.
[Schm2, p.~193]) guarantees that $\tilde S$ is absolutely irreducible. We
now apply
the Lang-Weil theorem, as stated e.g. in [Schm2, Thm.\ 5A, p.\ 210 or Cor.
5C, p.\ 213]
or [Se, p.\ 184]. We deduce that the number $N(p)$ of solutions $({\bf
x},y):=(x_0,\ldots ,x_r,y)\in\F_p^{r+2}$ to $\tilde S({\bf x},y)=0$
satisfies
$$
N(p)=p^{r+1}+O(p^{r+{1\over 2}}).
$$
Recall that  $d|D$, so $p\equiv 1\pmod d$. Hence, for each ${\bf
x}\in\F_p^{r+1}$
which corresponds to a solution $({\bf x},y)$ with $F({\bf x})\not\equiv
0\pmod{\cal
P}$, we find $d$ distinct solutions to be counted in $N(p)$. But the number of
solutions of  $F({\bf x})\equiv 0\pmod{\cal P}$ in $\F_p^{r+1}$ is trivially
$O(p^r)$. Therefore the number of ${\bf x}\in\F_p^{r+1}$ such that $\tilde
F({\bf
x})$ is a $d^{\rm th}$ power in $\F_p$ is  at most $(p^{r+1}/d)+O(p^{r+{1\over
2}})\le
2p^{r+1}/3$, say, for large enough $p$.\footnote{$^1$}{This argument goes back to M.
Eichler
and has been used by S.\ D. Cohen in the context of large sieve applied to
Hilbert's
Irreducibility Theorem. See [Se, Ch. 13].}

In particular, if $p$ has been chosen large enough, we may find
$${\bf x}=(x_0,\ldots ,x_r)\in\F_p^{r+1}$$ such that $x_0\cdots x_r\not =0$ and
$\tilde F(x_0,\ldots ,x_r)$ is not a $d^{\rm th}$ power in $\F_p$.  
\medbreak

Observe that $p$ splits completely in $\Q(\zeta_{p-1})$, so it splits
completely in
$k_3:=k_1(\zeta_{p-1})=k(\zeta_{p-1})$. Let $\cal Q$ be a prime ideal of $k_3$
above $\cal P$. Observe that the reduction of $\zeta_{p-1}$ generates the
group of
nonzero classes modulo $\cal Q$. So we may choose integers $h_0,\ldots
,h_r$ such
that $\zeta_{p-1}^{h_i}\equiv x_i\pmod{\cal Q}$. We determine $\zeta_m$ by
putting $\zeta_m=\zeta_{p-1}^D$. What we have proved implies the
following

\demo{{C}laim}
$\Psi := F(\zeta_{p-1}^{h_0},\ldots ,\zeta_{p-1}^{h_r})=
f(\zeta_{p-1}^{h_0},\zeta_m^{h_1},\ldots ,\zeta_m^{h_r})$
is not a $d^{\rm th}$ power in $k_3$.\enddemo

Consider now the field $K:=k(\zeta_{p-1},\beta^{1/m},\delta_1^{1/m},\ldots
,\delta_r^{1/m})$, where the $\delta_i$ are defined by (3). Plainly
$\Psi\in K$.

Suppose that $K$ has degree $<m$ over $k(\zeta_{p-1},\delta_1^{1/m},\ldots
,\delta_r^{1/m})$. Then,  Kummer theory for $K/k_3$ shows that $\beta$
belongs to the
multiplicative group generated by $\delta_1^{\ell\over m},\ldots
,\delta_r^{\ell\over m}$ and the $\ell^{\rm th}$ powers in $k(\zeta_{p-1})^*$, for
some prime $\ell$ dividing $m$. However this contradicts Lemma 2 (with
$M={m\over\ell}$), since each prime divisor of $m$ is $>L$ by property (a)
above.

Therefore $K/k(\zeta_{p-1},\delta_1^{1/m},\ldots ,\delta_r^{1/m})$ is cyclic of
degree $m$. We let $\tau$ be a generator of the Galois group.  

Recall that $\Psi$ is not a $d^{\rm th}$ power in $k_3= k(\zeta_{p-1})$. We contend
now that $\Psi$ is not a $d^{\rm th}$ power in $K$. In fact, the degree of $K$ over
$k(\zeta_{p-1})$ divides a power of $m$, while $k(\zeta_{p-1})(\Psi^{1/d})$ $k(\zeta_{p-1})$ has
degree $d$. But each prime factor of $m$ is $>L>d$, proving the
contention. 

Let $\xi$ be an element of ${\rm Gal}(K(\Psi^{1/d})/k(\zeta_{p-1}))$ which
induces $\tau$
on $K$ and is not the identity on $\Psi^{1/d}$. We let $\Pi$ be the set of
primes of
$k(\zeta_{p-1})$ of degree $1$ over $\Q$ and whose Frobenius in
$K(\Psi^{1/d})$ is
$\xi$. (Note that we are dealing with an abelian extension.) Cebotarev's
theorem guarantees that $\Pi$ is infinite. For $\pi\in\Pi$, let $\F_q$ be the
residue field of $\pi$. For almost all $\pi\in\Pi$ we may reduce
$\beta^{\pm 1}$,
$\gamma_i^{\pm 1}$, $\delta_i^{\pm 1}$, $f$ and $\Psi$ modulo $\pi$. Observe
that\medbreak

\item{(i)} Since the Frobenius of $\pi$ does not fix $\Psi^{1/d}$, we have that
$\Psi$
is not a $d^{\rm th}$ power modulo $\pi$.\smallbreak
\item{(ii)} Since $\xi$ induces $\tau$, which fixes $\delta_i^{1/m}$,  we have that
each $\delta_i$ is an $m^{\rm th}$ power modulo $\pi$, so

\vglue-16pt
$$
\delta_i^{q-1\over m}\equiv 1\pmod\pi .\eqno(5)
$$

\item{(iii)} Since $\xi$ induces $\tau$, which sends $\beta^{1/m}$ to
$\zeta_m^a\beta^{1/m}$ for some $a$ coprime with $m$, we have that

\vglue-16pt
$$
\beta^{q-1\over m}\equiv \zeta_m^a\pmod\pi .\eqno(6)
$$

\item{(iv)} Since $\pi$ has degree $1$ over $\Q$, there exists an integer
$n_0\in\Z$ such that $n_0\equiv\zeta_{p-1}\pmod\pi$.\medbreak

We can now obtain the desired contradiction as follows. By (5) and (6) we have
$$
\gamma_i^{q-1\over m}\equiv \beta^{h_i{q-1\over m}}\equiv
\zeta_m^{ah_i}\pmod\pi.
$$
Therefore, for an integer $a'$ such that $aa'\equiv 1\pmod m$,
$$
\gamma_i^{a'{q-1\over m}}\equiv \zeta_m^{h_i}\pmod\pi .
$$

Choose now an integer $n\in\N$ such that
$$
n\equiv n_0^{h_0}\pmod q ,\qquad n\equiv a'{q-1\over m}\pmod {q-1}.
$$
These
congruences imply
$$
\eqalign{
f(\zeta_{p-1}^{h_0},\zeta_m^{h_1},\ldots ,\zeta_m^{h_r})&\equiv
f(n_0^{h_0},\gamma_1^{a'{q-1\over m}},\ldots ,\gamma_r^{a'{q-1\over
m}})\cr
&\equiv f(n,\gamma_1^n,\ldots ,\gamma_r^n)=b(n)\pmod\pi .\cr}
$$

On the other hand the left side is $\Psi$ and the right side  is a $d^{\rm th}$ power
in $k$ by assumption. Since $b(n)$ is $\pi$-integral for almost all $\pi$, we
may reduce and get a contradiction with (i) above.

\demo{Added in proof} The principle of replacing powers $\alpha^{n}$ with roots of unity by reduction is used also in
[BBS].
\enddemo
 
\AuthorRefNames [Schm2]

\references

[BBS]
\name{D. Barsky, J.-P. B\`ezivin}, and \name{A. Schinzel}, Une caract\'erisation arithm\'etique de suites r\'ecurrentes
lin\'eaires, {\it J. reine angew.\ Math}.\ {\bf 494} (1998), 73--84.

[B]  \name{J.-P. B\'ezivin}, Sur les propri\'et\'es arithm\'etiques du produit de
Hadamard, {\it Ann.\ Fac.\ Sci.\ Toulouse Math.}\  {\bf 11} (1990), 7--19.

[BMZ]  \name{E. Bombieri, D. Masser}, and \name{U. Zannier},
Intersecting a curve with algebraic subgroups of multiplicative groups, 
{\it Int.\ Math.\ Res.\ Notices}, to appear.

[CZ]  \name{P. Corvaja} and \name{U. Zannier}, Diophantine equations with power sums and
universal Hilbert sets, {\it Indag.\ Mathem}.\ {\bf 9}  (1998),
317--332. 

[PZ] 
\name{A. Perelli} and  \name{U. Zannier}, Arithmetic properties of certain recurrent
sequences, {\it J. Austral.\ Math.\ Soc.}\ (Ser. A) {\bf 37} (1984),
4--16.

[P]  
\name{Y. Pourchet}, Solution du probl\`eme arithm\'etique du quotient de
Hadamard de deux fractions rationelles, {\it C. R. Acad.\ Sci.\ Paris}, 
S\'er.\ $A$ {\bf 288} (1979), A1055--A1057.

[RvdP] 
\name{R. Rumely} and \name{A.\ J. Van der Poorten}, A note on the Hadamard $k^{\rm th}$
root of a rational function, {\it J. Austral.\ Math.\ Soc.}\
(Ser. A) {\bf 43} (1987), 314--327.

[Schm1]  \name{W.\ M. Schmidt}, {\it Diophantine Approximation{\rm ,}
Lecture Notes in Math}.\ {\bf 785}, Springer-Verlag, New York, 1980.

[Schm2]  
\bibline, {\it Equations over Finite Fields. An
Elementary Appraoch{\rm , }
Lecture Notes in Math}.\ {\bf 536}, Springer-Verlag, New York, 1976.

[Se]  \name{J-P. Serre}, {\it Lectures on the Mordell-Weil Theorem{\rm , }
Aspects of Math}. {\bf E15}, Friedr.\ Viewer and  Sohn, Braunschweig,
1990.

[vdP1]  \name{A.\ J. Van der Poorten}, Some facts that should be better known,
especially  about rational functions, in {\it Number Theory  and
Applications}  (Banff, AB, 1988), 497--528,
Kluwer Acad.\ Publ., Dordrecht, 1989. 

[vdP2]  \bibline, Solution de la conjecture de Pisot sur le
quotient de Hadamard de deux fractions rationelles, {\it C. R. Acad.\ Sci.\
Paris},  S\'erie I {\bf 306} (1988), 97--102.

[vdP3] \bibline, A note on Hadamard roots of rational
functions, {\it Rocky Mountain J.\ Math}.\ {\bf 26} (1996),
1183--1197. 

\endreferences

\bye